\documentclass[12pt]{article}  
\usepackage{graphicx,amsmath,amsthm,amssymb,dsfont, mathrsfs }  
\numberwithin{equation}{section}
\def\ZZ{\mathbb Z}

\def\cH{\mathcal H}
\def\sL{\mathscr L}

\def\cF{\mathcal F}

\def\fa{\mathfrak a}
\def\fm{\mathfrak m}

\def\cP{\mathcal P}
\def\bcP{\boldsymbol{\cP}}

\def\bx{{\bf x}}
\def\by{{\bf y}}
\def\bz{{\bf z}}

\def\bm{{\bf m}}
\def\bk{{\bf k}}

\def\bu{{\bf u}}
\def\bw{{\bf w}}

\def\br{{\bf r}}

\def\btau{\boldsymbol{\tau}}

\def\bS{\boldsymbol{\Sigma}}
\def\bs{{\bf 0}}
\def\b1{{\bf 1}}
\def\d1{\mathds{ 1}} 
\def\mod{{\rm mod}}
\def\ad{{\rm and}}
\def\with{{\rm with}}
\def\where{{\rm where}}

\def\qed{ \ \vrule width.2cm height.2cm depth0cm\smallskip}

\begin{document}



\title{On the lower bound of the discrepancy  of Halton's sequence II}
\author{Mordechay B. Levin}

\date{}

\maketitle

\begin{abstract}
Let $ (H_s(n))_{n \geq 1} $ be an $s-$dimensional generalized Halton's sequence. Let $\emph{D}^{*}_N$ be the discrepancy of the sequence $ (H_s(n) )_{n = 1}^{N} $.
It is known that $D^{*}_{N} =O(\ln^s N)$ as $N \to  \infty $. 
 In this paper, we prove that this estimate  is exact. Namely, there exists a constant $C(H_s)>0$, such that
 $$
         \max_{1 \leq M \leq N}  
		M \emph{D}^{*}_{M} \geq C(H_s) \log_2^s N \quad {\rm for} \; \; N=2,3,... \; .
 $$
\end{abstract}
%
%
\section{Introduction }
%

 Let $(\beta_{n})_{n \geq 1}$ be a  sequence in  the unit cube $[0,1)^s$,
$B_{\by}=[0,y_1) \times \cdots \times [0,y_s) $,
\begin{equation}\label{I1}
\Delta(B_{\by}, (\beta_{n})_{n=1}^{N}  )= \sum\nolimits_{n=1}^{N}  ( \b1_{B_{\by}}(\beta_{n}) -  y_1 \cdots y_s), 
\end{equation}
where $\b1_{B}(\bx) =1, \; {\rm if} \;\bx  \in B$, 
and $   \b1_{B_{\by}}(\bx) =0,$  if $ 
\bx   \notin B$.

We define the star {\it discrepancy} of an 
$N$-point set $(\beta_{n})_{n=1}^{N}$ as
\begin{equation} \label{2}
   \emph{D}^{*}((\beta_{n})_{n=1}^{N}) = 
    \sup\nolimits_{ 0<y_1, \ldots , y_s \leq 1} \; | 
  \Delta(B_{\by},(\beta_{n})_{n=1}^{N})/N |.
\end{equation}

In 1954, Roth proved that 
\begin{equation} \nonumber
   \limsup\nolimits_{N \to \infty } N (\ln N)^{-\frac{s}{ 2}} \emph{D}^{*}((\beta_{n})_{n=1}^{N})>0 . 
\end{equation}
According to the well-known conjecture (see, e.g., [BeCh, p.283]), this estimate can be improved
 to
\begin{equation}   \label{4}
 \limsup\nolimits_{N \to \infty } N (\ln N)^{-s} \emph{D}^{*}((\beta_{n})_{n=1}^{N})>0 . 
\end{equation}
In 1972, W. Schmidt proved this conjecture for $ s=1 $. For $s=2$, Faure and Chaix [FaCh] proved (\ref{4}) for a class of $(t,s)-$sequences.
See [Bi] for the most important results on this conjecture. 

 \texttt{ Definition.} {\it  An $s$-dimensional sequence  $((\beta_{n})_{n \geq 1})$ is of 
 \texttt{low discre-\\pancy} (abbreviated
l.d.s.) if $ \emph{D}^{*}((\beta_{n})_{n=1}^{N})=O(N^{-1}(\ln
N)^{s}) $ for $ N \rightarrow \infty $. }

Let $ p\geq 2 $ be an integer
 \begin{equation}\nonumber
 n=\sum_{j\geq 1}e_{p,j}(n) p^{j-1},\;  e_{p,j}(n) \in \{0,1, \ldots
 ,p-1\}, \;  {\rm and}    \; \phi_p(n)= \sum_{j\geq 1}e_{p,j}(n) p^{-j}.
 \end{equation}
Van der Corput    proved that $ (\phi_p(n))_{n\geq 0}$ is a $1-$dimensional l.d.s. (see [VC]). Let
\begin{equation}\nonumber
  \hat{H}_s(n)= (\phi_{\hat{p}_1}(n),\ldots ,\phi_{\hat{p}_s}(n)), \quad n=0,1,2,...,
\end{equation}
 where $
\hat{p}_1,\ldots ,\hat{p}_s\geq 2 $ are pairwise coprime integers.
Halton   proved that $ ( \hat{H}_s(n))_{n\geq 0}$ is an $s-$dimensional l.d.s. (see [Ha]). 
For other examples of  l.d.s. see e.g. in [BeCh], [FKP], [Ni].
In [Le2] we proved that Halton's sequence  satisfies  (\ref{4}). In this paper we generalize this result.

Let $Q=(q_1,q_2,....)$ and $Q_j=q_1q_2....q_j$, where $q_j\ge2$ \ $(j=1,2,\dots)$ is a sequence of integers.
Consider Cantor's expansion of $x \in [0,1):$
$$x=\sum\nolimits_{j=1}^{\infty} x_j / Q_j, \quad
x_j \in \{0,1,\dots,q_j-1\}, \quad x_j \neq q_j-1\; {\rm for \;infinitely \; many} \;j.$$
The $Q-$adic representation of $x$ is then unique.
We define the odometer transform 
\begin{equation} \label{T} 
    T_Q(x) := (x_k+1)/ Q_k + \sum\nolimits_{j \geq k+1}  x_j / Q_j, \qquad \;\;
		T_Q^n(x)=T_Q(T_Q^{n-1}(x)),
\end{equation}
$n=2,3,\dots,T_Q^0(x)=x,$ where $k= \min \{j \;| \; x_j \ne q_i-1 \}$.

For $Q=(q,q,....)$, we obtain  von Neumann-Kakutani's $q-$adic 
adding machine (see, e.g., [FKP]). 
As is known, the  sequence $(T_{Q}^n(x))_{n\ge1}$ coincides for $x=0$  with the van
der Corput sequence (see e.g., [FKP, \S 2.5]).

Let  $h_i \geq 1$, $q_{i,j} \geq 2$ be integers $(1 \leq j \leq h_i, 1 \leq i \leq s)$, $p_{i,j} \in \{q_{i,1},..., q_{i,h_i} \}$, $g.c.d.(q_{i,k},q_{j,l})=1$ for $i\neq j$, 
 $\cP_i = (p_{i,1},p_{i,2},...)$,  $\bcP=(\cP_1,...,\cP_s)$,
\begin{equation} \nonumber
 \tilde{P}_{i,0} =1, \;\;  \tilde{P}_{i,j} =\prod_{ 1 \leq k \leq j}  p_{i,k}, \; i \in [1,s],\; j \geq 1,
 \quad 
T_{\bcP}(\bx) = (T_{\cP_1}(x_1),...,T_{\cP_s}(x_s)),
\end{equation}
\begin{equation}   \label{1.4}
n=\sum_{j \geq 1} e_{p_{i,j},j}(n) \tilde{P}_{i,j-1},
 \quad e_{p_{i,j},j}(n) \in \{0,1, \ldots
 ,p_{i,j} -1\},\; n=0,1,... \;,
 \end{equation}
\begin{equation} \label{1.4a} 
 \varphi_{\cP_i}(n)=     \sum_{j \geq 1}e_{p_{i,j},j}(n) \tilde{P}^{-1}_{i,j}, \;\; \quad  H_{\bcP}(n)= (\varphi_{\cP_1}(n),\ldots ,\varphi_{\cP_s}(n)).
\end{equation}
We note that $H_{\bcP}(n) =T^n_{\bcP}(\bs)$ for $n=0,1,...$ .
 
Let $\Sigma_i =(\sigma_{i,j})_{j \geq 1}$ be a sequence  of corresponding permutations 
$\sigma_{i,j}$ of $\{0,1,...,p_{i,j}-1\}$ for $j \geq 1$, 
 $\bS=(\Sigma_1,...,\Sigma_s)$,  $\bx=(x_1,...,x_s)$,
\begin{equation} \nonumber
  \tilde{\bS}(\bx) =(\tilde{\Sigma}_1(x_1),...,\tilde{\Sigma}_s(x_s)), \quad  \tilde{\Sigma}_i(x_i) = \sum_{j \geq 1} \sigma_{i,j}(x_{i,j})/\tilde{P}_{i,j},\; 
	x_i=\sum_{j \geq 1} x_{i,j}/\tilde{P}_{i,j}.
\end{equation}
We consider the following generalization of the  Halton sequence (see [Fa],[He],  [FKP]):
\begin{equation} \label{1.4d} 
  H_{\bcP}^{\bS}(n,\bx)= \tilde{\bS}(T_{\bcP}^n(\bx)), \quad  \; n=0,1,2,... \; .
\end{equation}
We note that $(H_{\bcP}^{\bS}(n,\bx))_{n \geq 0} $ coincides for $\bx=\bs$ and $s=1$ with the 
Faure sequence $S_Q^{\Sigma}$ [Fa].
Similarly to [Ni, p.29-31], we get that $ (H_{\bcP}^{\bS}(n,\bx))_{n \geq 0}$ is of low 
discrepancy.

\section{The Theorem and its proof}
In this section we will prove \\

{\bf Theorem.} {\it Let  $s \geq 2$,   $h_0 = \max_i h_i$,  $q_0 = \max_{i,j} q_{i,j}$,
 $C_1=2s h_0 q_0^{s} \log_2 q_0$ and  $C =2^{s+3}s^s h_0^s q_0^{s^2} \log_2^{s} q_0$,
 $\log_2 N \geq 2q_0^{s-1} C_1$. Then }
\begin{equation}\label{Th}
    \inf_{\bx\in [0,1)^s}  \max_{1 \leq M \leq  N} 
			M \emph{D}^{*}((  H_{\bcP}^{\bS}(n,\bx) )_{n=1}^{M}) \geq C^{-1} \log_2^{s } N   .
\end{equation}   
This result  supports  conjecture (\ref{4}) (see also  [Le1] and [Le3]).

First we will construct a double sequence $(\tau_{i,j})_{1 \leq i \leq s, j \geq 1}$. In order to construct $(\tau_{i,j})_{1 \leq i \leq s, j \leq 1}$, we define auxiliary sequences 
$\sL_{i,j}^{(\fm)}, L^{(\fm)}_i ,l_{i,j}, \cF_{i,b}^{(\fm)},...$
  as follows.\\

 {\sf 2.1  Construction of the sequence $(\tau_{i,j})$.}\\
Let  $\fm=[ \log_{q_0}(N)/s -1]$ with $q_0 = \max_{i,j} q_{i,j}$, 
$a_{i,j} \equiv  \sigma^{-1}_{i,j}(0) - \sigma^{-1}_{i,j}(1) \; (\mod \; p_{i,j})$,
$a_{i,j} \in \{1,...,p_{i,j}-1\} $, 
\begin{eqnarray} \label{2.0}
     \sL_{i,j, \tilde{\fa}_{i,j}}^{(\fm)}  = \{ 1 \leq k \leq \fm \; | \; p_{i,k} =q_{i,j}, \; 
		a_{i,k} =\tilde{\fa}_{i,j}      \},   \\
	L^{(\fm)}_i =\#\sL_{i,g_{i,\fm}, \fa_{i,m} }^{(\fm)} = 
	\max_{1 \leq j \leq h_i, 1 \leq \tilde{\fa}_{i,j} < q_{i,j}} \# \sL_{i,j, \tilde{\fa}_{i,j}}^{(\fm)}, \quad
	\where \; g_{i,\fm} \in [1,h_i],   \nonumber
\end{eqnarray}  
$\fa_{i} =\fa_{i,m} \in [1, q_{i,g_{i,\fm}  }-1] $, $1 \leq i \leq s$.
We enumerate the set $\sL_{i,g_{i,\fm}, \fa_{i}}^{(\fm)}$:
\begin{equation}  \nonumber
 \sL_{i,g_{i,\fm, \fa_{i}}}^{(\fm)} = \{ l_{i,1} < \cdots <  l_{i, L^{(\fm)}_i} \}.
\end{equation}
 We see that 
\begin{equation} \label{2.4}
   L^{(\fm)}_i  \geq \fm/(h_iq_0) \qquad \ad \qquad a_{i, l_{i,j}} =\fa_{i}, \; i=1,...,s,\; j =1,...,\fm.
\end{equation}

Let $p_i = p^{(\fm)}_i = q_{i,g_{i, \fm} }$, $p_0 =p^{(\fm)}_0= p_1p_2 \cdots p_s$,  $\dot{p}_i =p_0/p_i$ and
\begin{equation}\label{2.5}
     \cF_{i,b}^{(\fm)}  = \{ 1 \leq k \leq L^{(\fm)}_i \; | \; \tilde{P}_{i,l_{i,k}}^{-1} \equiv  b \;\; ( {\rm mod} \;   \dot{p}_i  )    \}. 
\end{equation} 
We define $ F_i, \; m$ and $b_i=b_i^{(\fm)}$  as follows:
\begin{equation}\label{2.2}
        F_i =F^{(\fm)}_i =\#\cF_{i,b_i}^{(\fm)} = \max_{0 \leq b < \dot{p}_i}   \# \cF_{i,b}^{(\fm)}, \qquad m =\min_{1 \leq i \leq s} F^{(\fm)}_i .
\end{equation} 
It is easy to see that
\begin{equation}\label{1.31}
   m  \geq \min_{1 \leq i \leq s}  \fm/(h_i q_0 \dot{p}_i)  \geq \fm h_0^{-1}  q_0^{-s}
	 \geq C_1^{-1} \log_2 N,
\end{equation}
with $C_1=2s h_0 q_0^{s} \log_2 q_0$. 
We enumerate the set $F_{i,b_i}^{(\fm)}$:
\begin{equation}\nonumber
 \cF_{i,b_i}^{(\fm)}= \{ f_{i,1} < \cdots <  f_{i, F_i} \} .
\end{equation}
Let $\bk =(k_1,...,k_s)$, $\tau_{i, j} =l_{i,f_{i,j}}$, $\btau_{\bk} =(\tau_{1, k_1},...,
 \tau_{s, k_s})$,    $P_{i,k} = \tilde{P}_{i,\tau_{i,k}}$,
\begin{equation} \label{1.20}
  P_{\bk} = \prod_{i=1}^{s} P_{i,k_i}, \;\;
  M_{i,\bk} =\tilde{M}_{i,\btau_{\bk}}, \;\; {\rm with} \;\;
\tilde{M}_{i,\bk} \equiv \prod_{1\leq j \leq s, j \neq i} \tilde{P}_{j,k_j}^{-1}   \; ({\rm mod} \; \tilde{P}_{i,k_i}).
\end{equation}
By (\ref{2.5}), we have that $(b_i,\dot{p}_i) =1$ and $(b_j,p_i) =1$ for $i\neq j $ $(i,j=1,...,s)$.\\ 
Let $ c_i \equiv \prod_{1 \leq j \leq s,j \neq i} b_j \; ({\rm mod} \; p_i).$ 
According to (\ref{2.4}), (\ref{2.5}) and (\ref{1.20}), we obtain
\begin{equation} \label{1.50}
 (c_i,p_i)=1, \quad
  M_{i,\bk} \equiv c_i   \;\; ({\rm mod} \; p_i) \;\; {\rm and} \;\;
	a_{i, \tau_{i, j}}=\fa_{i}, \; j \geq 1, \; i=1,...,s.
\end{equation}
Let 
\begin{equation} \nonumber
   \tilde{p}_i=g.c.d.(\fa_{i},p_i), \quad \;\; \hat{p}_i=p_i/\tilde{p}_i, \quad \;\;
\hat{a}_i =\fa_{i}/\tilde{p}_i, \quad \;\; d_i \equiv c_i \fa_{i} \; (\mod \;\hat{p}_i),
\end{equation}
 $d_i \in \{1,..., \hat{p}_i-1 \}$. Hence
\begin{equation} \label{1.50a}
   d_i/ \hat{p}_i \equiv c_i \fa_{i}/p_i \; (\mod \;1), \quad ( d_i, \hat{p}_i)=1, \; \ad \;
	\hat{p}_i >1,  
	 \; i=1,...,s.
\end{equation}
Let  $\bm =(m,...,m)$.   From (\ref{1.4}) and (\ref{1.20}), we derive
\begin{equation} \label{1.60}
  2P_{\bm }  \leq   2\prod\nolimits_{i=1}^s \prod\nolimits_{j=1}^{\tau_{i, m}}  p_{i,j} \leq 2q_0^{\fm s}   
	 \leq q_0^{s[s^{-1} \log_{q_0} N]} \leq N.
\end{equation}\\


 {\sf 2.2 Using the Chinese Remainder Theorem.}\\
Let $x_i =\sum_{j \geq 1}  x_{i,j}\tilde{P}_{i,j}^{-1}$, with $x_{i,j} \in \{0,1,...,p_{i,j}-1  \}$, 
 $i=1,...,s$.  
We define the truncation 
\begin{equation}  \nonumber
        [x_i]_r =\sum_{1 \leq j \leq r}  x_{i,j}\tilde{P}_{i,j}^{-1} \quad \with \quad r \geq 1.
\end{equation}
If $x = (x_1, . . . , x_s)  \in [0, 1)^s$, then the truncation $[\bx]_{\br}$ is defined coordinatewise, that is, $[\bx]_{\br}= 
( [x_1]_{r_1}, . . . , [x_s]_{r_s})$, where $\br =(r_1,...,r_s)$.

By (\ref{1.4a}), we have  
\begin{equation}\nonumber
    [\varphi_{\cP_i}(n)]_{r_i}=[x_i]_{r_i} \;
			\Leftrightarrow  \; n \equiv    \sum_{1 \leq j \leq r}  x_{i,j}\tilde{P}_{i,j-1}  \;\; ({\rm mod} \; \tilde{P}_{i,r}).
\end{equation}
Applying (\ref{1.20}) and the Chinese Remainder Theorem, we get
\begin{equation} \label{2.16a}
  [H_{\bcP}(n)]_{\br} = [\bx]_{\br}  \; \Longleftrightarrow  \; n \equiv   \check{x}_{\br}  \; ({\rm mod} \; \tilde{P}_{\br}),
\end{equation}
\begin{equation}\label{2.16}
 \check{x}_{\br} \equiv \sum_{i=1}^s  \tilde{M}_{i,\br}  
   \tilde{P}_{\br}\tilde{P}_{i,r_i}^{-1} 
  \sum_{1 \leq j \leq r}  x_{i,j}\tilde{P}_{i,j-1} \; (\mod \; \tilde{P}_{\br} ), \; \qquad 
	\check{x}_{\br} \in [0,\tilde{P}_{\br}).
\end{equation}
It is easy to verify that if $r_i^{'} \geq r_i$, $i=1,...,s$, then
\begin{equation} 
              \check{x}_{\br^{'} } \equiv \check{x}_{\br}  \;\; ({\rm mod} \;  \tilde{P}_{\br} ).
\end{equation}
According to (\ref{T}), we get 
\begin{equation}\nonumber
{\rm if} \; [\bw]_{\br}= [\bx]_{\br}, \quad {\rm then} \quad [T_{\bcP}^n(\bw)]_{\br}= [T_{\bcP}^n(\bx)]_{\br}, \qquad   n=0,1,... \;.
\end{equation}
From  (\ref{T}), (\ref{1.4a}) and (\ref{2.16a}), we obtain
\begin{equation}\nonumber 
   [T_{\bcP}^W(\bs)]_{\br}= [H_{\bcP}(W)]_{\br} = [\bx]_{\br},  \; \qquad \; W=\check{x}_{\br}.
\end{equation}
Hence
\begin{equation}\nonumber 
[T_{\bcP}^n(\bx)]_{\br} =[T_{\bcP}^n(T_{\bcP}^W(\bs))]_{\br} =[T_{\bcP}^{n+W}(\bs)]_{\br}= 
        [H_s(n+W)]_{\br} .
\end{equation}
Let
\begin{equation} 
     W_{\bm}(\bx):=          \check{x}_{\bm} \in [0,P_{\bm}) .
\end{equation}
Therefore
\begin{equation} \label{1.15}
[T_{\bcP}^n(\bx)]_{\br} = [H_{\bcP}(n+W_{\bm}(\bx))]_{\br} \;\;  \quad \;\;  1 \leq r_i \leq m, 
 \; 1 \leq i \leq s, \; n \geq 0 .
\end{equation}\\

 {\sf 2.3  Construction of boundary points $y_1,...,y_s$ and $u_1,...,u_s$.}\\
Let $\by =(y_1,...,y_s)$ with $y_i = \sum_{1 \leq j \leq m} P_{i,j}^{-1}$, and let   $\ddot{y}_{i,k_i}=\sum_{1 \leq j \leq k_i} P_{i,j}^{-1} $, $k_i \geq 1$,  $i=1,...,s$, 
 $\bk =(k_1,...,k_s)$, 
\begin{equation} \label{1.15a}
B_{\by}=[0,y_1) \times \cdots \times [0,y_s)  , \quad B^{(\bk)} =
\prod\nolimits_{i=1}^s \; [\ddot{y}_{i,k_i} -P_{i,k_i}^{-1},\ddot{y}_{i,k_i}).
\end{equation}
We deduce
\begin{equation}\label{1.16}
    B_{\by} =\bigcup_{k_1,...,k_s=1}^m B^{(\bk)}, 		\; {\rm and} \; 
		\b1_{B_{\by}} (\bz) -  y_1 \cdots y_s = 		
		\sum_{k_1,...,k_s=1}^m    ( \b1_{B^{(\bk)}} (\bz) -  P_{ \bk}^{-1}).
\end{equation} 
Let $\bu =(u_1,...,u_s)$,  $u_i =\sum_{j \geq 1}^{\tau_{i,m}}  u_{i,j}\tilde{P}_{i,j}^{-1}$  
 with $u_{i,j} =\sigma_{i,j}^{-1}(y_{i,j}) $, $u_{i,j}^{*} =\sigma_{i,j}^{-1}(0)$, 
\begin{equation}\label{1.18}
   \bu^{(\bk)} =(u_1^{(k_1)},...,u_s^{(k_s)}) \;\;\; \with \;\;\; u_i^{(k_i)} =
	\sum_{j = 1}^{\tau_{i,k_i}-1}
		u_{i,j}\tilde{P}_{i,j}^{-1} + u_{i,\tau_{i,k_i}}^{*} \tilde{P}_{\tau_{i,k_i}}^{-1}, 
\end{equation}
\begin{equation} \nonumber
     \check{u}^{(\bk)}
		\equiv \sum_{i=1}^s  M_{i,\bk}  
    P_{\bk} P_{i,k_i}^{-1} 
  \Big( \sum_{j=1}^{\tau_{i,k_i}-1}  u_{i,j}\tilde{P}_{i,j-1} 
	  + u_{i,{\tau_{i,k_i}}}^{*} \tilde{P}_{i,{\tau_{i,k_i}-1}} \Big) \; (\mod \; P_{\bk} ), 
\end{equation}
\begin{equation} \nonumber
     \check{u}_{\bk}
		\equiv \sum_{i=1}^s  M_{i,\bk}  
   P_{\bk} P_{i,k_i}^{-1} 
   \sum_{j=1}^{\tau_{i,k_i}}   u_{i,j}\tilde{P}_{i,j-1} 
	  \; (\mod \; P_{\bk} ), \;\;   \quad  
	\check{u}^{(\bk)}, \check{u}_{\bk} \in [0,P_{\bk}).
\end{equation}

According to (\ref{2.0})-(\ref{1.20}), we have $ p_{i, \tau_{i,k_i}}= p_i$,
 $k_i=1,...,m$, $i=1,...,s$. \\
By (\ref{2.0}), we get $a_{i,\tau_{i,k_i}} \equiv  \sigma^{-1}_{i,\tau_{i,k_i}}(0) - \sigma^{-1}_{i,\tau_{i,k_i}}(1) 
 \equiv   u_{i,\tau_{i,k_i}}^{*} - u_{i,\tau_{i,k_i}} 
\;({\rm mod} \; p_i)$. \\
From (\ref{1.50}), we obtain $ a_{i, \tau_{i,k_i}}= \fa_i$,
 $k_i=1,...,m$, $i=1,...,s$. 
Hence
\begin{equation}\label{1.19a}
       \check{u}^{(\bk)} \equiv  \check{u}_{\bk} +A_{\bk}\;
		 (\mod \; P_{\bk} ), \; \where  \;  A_{\bk} \equiv 
				\sum\nolimits_{i=1}^s  M_{i, \bk}   P_{ \bk} p_i^{-1} \fa_i \; ({\rm mod} \; P_{\bk})   
\end{equation}
with $A_{\bk} \in[0,  P_{ \bk})$. 

Let $\bw=(w_1,...,w_s):= H_{\bcP}^{\bS}(n,\bx)=\tilde{\bS}(T_{\bcP}^n(\bx))$. \\
We see from (\ref{1.15a}) and (\ref{1.18}) that
\begin{equation}\nonumber
			\bw  \in B^{(\bk)}
		\Leftrightarrow w_{i,j}  =y_{i,j},\;  j \in [1, \tau_{i,k_i}),\;
		w_{i,\tau_{i,k_i}}  =0, \;		i\in [1,s] 
		\Leftrightarrow \sigma_{i,j}(w_{i,j})  =u_{i,j}
\end{equation}
\begin{equation}\nonumber
 	 1 \leq j \leq \tau_{i,k_i} -1,\; 
		\sigma_{i,j}(w_{i,\tau_{i,k_i}}) =u_{i,\tau_{i,k_i}}^{*} , i=1,...,s \;
		\Leftrightarrow [T_{\bcP}^n(\bx)]_{\btau_{\bk}} = \bu^{(\bk)}.	
\end{equation}
Applying  (\ref{2.16a}), (\ref{2.16}), (\ref{1.15}), (\ref{1.18}) and (\ref{1.19a}), we have  
\begin{equation}\nonumber
      H_{\bcP}^{\bS}(n,\bx) \in B^{(\bk)} \Leftrightarrow [T_{\bcP}^n(\bx)]_{\btau_{\bk}} = \bu^{(\bk)} \Leftrightarrow [H_{\bcP}(n+W_{\bm}(\bx))]_{\btau_{\bk}} = \bu^{(\bk)} 
\end{equation}
\begin{equation}\nonumber
  \Leftrightarrow 	n+W_{\bm}(\bx) \equiv   \check{u}^{(\bk)}
			\;\; ({\rm mod} \; P_{\bk})      \Leftrightarrow  n \equiv v_m +A_{\bk} \;\; ({\rm mod} \; P_{\bk}),  
\end{equation}
where $v_m \equiv -W_{\bm}(\bx)  + \check{\bu}_{\bm}
\equiv -W_{\bm}(\bx)  + \check{\bu}_{\bk} \; ({\rm mod} \; P_{\bk})$ and 
$v_m \in [0,P_{\bm})$.\\
Hence
\begin{equation}  \label{1.70}
     H_{\bcP}^{\bS}(n,\bx) \in B^{(\bk)}     \Longleftrightarrow   n \equiv   v_{m} +A_{\bk}
      \; ({\rm mod} \; P_{ \bk}), \quad  v_m \in [0,P_{\bm}), \; n \geq 0.  
\end{equation}\\ 
{\bf Completion of the proof of Theorem.} 

{\bf Lemma 1.} {\it Let
\begin{equation} \label{1.32b}
\alpha_{m}:=  \frac{1}{P_{ \bm}}  \sum_{ M=1}^{ P_{\bm}} 
  \Delta(B_{\by}, ( H_{\bcP}^{\bS}(n,\bx))_{n=v_m}^{v_m+M-1} )  .
\end{equation}
Then 
\begin{equation}  \label{1.23a}
\alpha_{m}  =  \sum_{1 \leq k_1,...,k_s \leq m}  \Big(\frac{1}{2} -   \frac{A_{\bk}}{P_{ \bk}}  - \frac{1}{2P_{ \bk}} \Big).
\end{equation}}

{\bf Proof.}
Let $\cH_n :=H_{\bcP}^{\bS}(n,\bx)$. Using (\ref{1.70}), we have

\begin{equation} \label{1.22}  
        \sum_{n =  v_m + M_1 P_{ \bk} }^{ v_m +(M_1+1)P_{ \bk}-1}  ( \b1_{B^{(\bk)}} (\cH_n) - P_{ \bk}^{-1} )  =0
\end{equation}
and
\begin{equation} \nonumber
				 \sum_{n= v_m + M_1 P_{ \bk}}^{ v_m +M_1P_{ \bk}+M_2-1}  ( \b1_{B^{(\bk)}} (\cH_n) - P_{ \bk}^{-1} )  =\sum_{ n \in [v_m, v_m +M_2)}  ( \b1_{B^{(\bk)}} (\cH_n) - P_{ \bk}^{-1} )
\end{equation}
\begin{equation} \nonumber
				= \sum_{  n \in [v_m, v_m +M_2), n= v_m +A_{\bk} }  1 - M_2P_{ \bk}^{-1}  
       =  \b1_{[0,M_2)}(A_{\bk}) -M_2P_{ \bk}^{-1},
\end{equation}
with $M_1 \geq 0$ and $M_2 \in [0, P_{\ \bk})$, $M_1,M_2 \in \ZZ$.

From (\ref{I1}) and (\ref{1.16}), we get
\begin{equation}\nonumber
        \Delta(B_{\by}, (\cH_n)_{n=v_m}^{v_m+M-1} ) 
  =      \sum_{ n=v_m}^{v_m +M-1}  ( \b1_{B_{\by}}(\cH_n) -  y_1 \cdots y_s)  
\end{equation}
\begin{equation}\label{1.32}
      = \sum_{ k_1,...,k_s=1}^{ m} \rho(\bk,M), \quad {\rm with} \quad  
 \rho(\bk,M)= 
       \sum_{ n=v_m}^{v_m +M-1}  ( \b1_{B^{(\bk)}} (\cH_n) -  P_{ \bk}^{-1})  .
\end{equation}

By (\ref{1.32b}), we obtain
\begin{equation} \label{1.32a}
\alpha_{m}   = \sum_{1 \leq k_1,...,k_s \leq m} \alpha_{m,\bk}, \;\; {\rm with} \;\;  
 \alpha_{m,\bk}= \frac{1}{P_{ \bm}} \sum_{ M=1}^{ P_{ \bm}} \rho(\bk,M).
\end{equation}
Bearing in mind  (\ref{1.22})-(\ref{1.32}), we derive
\begin{equation}\nonumber
  \alpha_{m,\bk} =
      \frac{1}{P_{ \bm}} \sum_{ M_1=0}^{ P_{ \bm}/ P_{ \bk} -1}\sum_{ M_2=1}^{ P_{ \bk}}  \Big(
       \sum_{ n=v_m}^{v_m  +M_1P_{ \bk} -1}  ( \b1_{B^{(\bk)}} (\cH_n) -  P_{ \bk}^{-1})  
\end{equation}
\begin{equation}\nonumber
    +   \sum_{ n=v_m  +M_1P_{ \bk}}^{v_m  +M_1P_{ \bk}+M_2-1}   ( \b1_{B^{(\bk)}} (\cH_n) -  P_{ \bk}^{-1})  \Big)
     =  
      \frac{1}{P_{\bm}} \sum_{ M_1=0}^{ P_{\bm}/ P_{ \bk} -1}\sum_{ M_2=1}^{ P_{ \bk}}  \Big(    \b1_{[0,M_2)}(A_{\bk}) -M_2P_{ \bk}^{-1}   \Big) 
\end{equation}
\begin{equation}\nonumber
        = \frac{1}{P_{ \bk}} \sum_{ M_2=1}^{ P_{ \bk}}  \Big(    \b1_{[0,M_2)}(A_{\bk}) -M_2P_{ \bk}^{-1}   \Big) = \frac{P_{ \bk} -A_{\bk}}{P_{ \bk}}-\frac{P_{ \bk}(P_{ \bk}+1)}{2P_{ \bk}^2} = \frac{1}{2} -   \frac{A_{\bk}}{P_{ \bk}}  - \frac{1}{2P_{ \bk}}.
\end{equation}
Using  (\ref{1.32a}), we have
\begin{equation}  \nonumber
\alpha_{m}  =  \sum_{1 \leq k_1,...,k_s \leq m}  \Big(\frac{1}{2} -   \frac{A_{\bk}}{P_{ \bk}}  - \frac{1}{2P_{ \bk}} \Big).
\end{equation}
Hence Lemma 1 is proved. \qed \\

{\bf Lemma 2.} {\it With notations as above,
\begin{equation}
|\alpha_{m}| 
 \geq \frac{m^s}{4p_0}
     \quad {\rm for} \quad m\geq2p_0.
\end{equation}}

{\bf Proof.}
From (\ref{1.50}) and (\ref{1.19a}), we get 
\begin{equation}\nonumber
  [0,1) \ni \frac{A_{\bk}}{P_{ \bk}} \equiv    \sum_{1 \leq i \leq s}  M_{i, \bk}   P_{ \bk} p_i^{-1}\fa_i/P_{ \bk} \; \equiv \;\frac{c_1 \fa_1}{p_1} + \cdots  + \frac{c_s\fa_s}{p_s}   \;\; ({\rm mod} \; 1). 
\end{equation}
Applying (\ref{1.50a}) and (\ref{1.23a}), we derive
\begin{equation}\label{1.23}
\alpha_{m}  
 = m^s \Big( \frac{1}{2} 
      - \{\alpha  \} \Big) - 
      \sum_{1 \leq k_1,...,k_s \leq m} \frac{1}{2P_{ \bk}},  \quad {\rm with} \quad  \alpha =  \frac{d_1}{\hat{p}_1} + \cdots  + \frac{d_s}{\hat{p}_s},
\end{equation}
where $( d_i, \hat{p}_i)=1, \; \hat{p}_i >1, 	 \; i=1,...,s.$ and $\{x\}$ is the fractional part of  $x$.
We have that if $\hat{p}_0 =\hat{p}_1\hat{p}_2\cdots \hat{p}_s \not \equiv 0 \;\; ({\rm mod} \; 2)$ then $ \alpha   \not \equiv 1/2 \;\; ({\rm mod} \; 1)$.
Let $\hat{p}_{\nu} \equiv 0 \;\; ({\rm mod} \; 2)$ for some $\nu \in [1,s]$, and let 
$ \alpha    \equiv 1/2 \;\; ({\rm mod} \; 1)$. Then 
\begin{equation} \nonumber
        (\hat{p}_{\nu}/2-d_{\nu})/p_{\nu}              \equiv \sum_{1 \leq i \leq s, \;i \neq \nu}  d_i/\hat{p}_i  \;\; ({\rm mod} \; 1)  \quad 
          {\rm and} \quad  a_1 \equiv a_2  \;\; ({\rm mod} \; p_0), 
\end{equation}
with $ a_1 = \hat{p}_0(\hat{p}_{\nu}/2-d_{\nu})/\hat{p}_{\nu} $  and $a_2 = \sum_{i \neq \nu}\hat{p}_0d_i/\hat{p}_i $. 
Let $j \in [1,s]$ and $j \neq \nu$. 
We see that $a_1  \equiv 0 \;\; ({\rm mod} \;\hat{p}_j)$ and $a_2  \not \equiv 0 \;\; ({\rm mod} \; \hat{p}_j)$.
We get a contradiction. Hence  $ \alpha   \not \equiv 1/2 \;\; ({\rm mod} \; 1)$.
We have
\begin{equation}\nonumber
 \Big| \frac{1}{2} - \Big\{  \alpha\Big\} \Big| =
     \Big| \frac{1}{2} - \Big\{   \Big(\frac{d_1}{\hat{p}_1} + \cdots  + \frac{d_s}{\hat{p}_s} \Big) \Big\} \Big| =  \frac{|a|}{2\hat{p}_0}, \quad {\rm with \;some \;integer} \; a. 
\end{equation}
Thus $  | 1/2 - \{   \alpha \} |  \geq 1/(2\hat{p}_0) \geq 1/(2 p_0)  $ with $p_0=p_1...p_s$, 
 $(p_0,\hat{p}_0)=\hat{p}_0$. 

Bearing in mind that $P_{ \bk} \geq 2^{k_1+k_2 + \cdots +k_s}  $, we obtain from (\ref{1.23}) that
\begin{equation}\label{1.30}
|\alpha_{m}| \geq   
  \frac{m^s}{2p_0}  -\frac{1}{2} =  \frac{m^s}{2p_0} (1 - \frac{p_0}{ m^{s}}) \geq \frac{m^s}{4p_0}
     \quad {\rm for} \quad m\geq2p_0 .
\end{equation}
Hence Lemma 2 is proved. \qed

Going back to the proof of Theorem, by (\ref{Th}) and (\ref{1.31}), we get 
\begin{equation}\nonumber
  m^s (4p_0)^{-1} \geq (4p_0)^{-1} C_1^{-s} \log_2^s N =2C^{-1}  \log_2^s N , \;\; \ad \;\;
	m \geq C_1^{-1} \log_2 N  \geq 2p_0,
\end{equation}
where $C_1=2s h_0 q_0^{s} \log_2 q_0$ and $C = (8p_0)^{-1} C_1^s
=2^{s+3}s^s h_0^s q_0^{s^2} \log_2^{s} q_0$.

Using  (\ref{1.60}) and (\ref{1.70}), we have that $v_m + P_{\btau m} \leq 2P_{\bm}   \leq N$.\\
According to (\ref{1.30}), (\ref{1.32b}) and (\ref{2}),  we obtain
\begin{equation}\nonumber
2C^{-1}\log_2^s N   \leq   m^s (4p_0)^{-1} \leq  |\alpha_{m} |   \leq  \sup_{1  \leq M \leq P_{\bm}}   M \emph{D}^{*}(( H_{\bcP}^{\bS}(n,\bx))_{n=v_m}^{v_m+M-1})     
\end{equation}  
\begin{equation}\nonumber
  \leq  \sup_{1 \leq L,L+M \leq 2P_{\bm}}   M \emph{D}^{*}(( H_{\bcP}^{\bS}(n,\bx))_{n=L}^{L+M-1})   \leq  2\sup_{1 \leq M \leq N}   M \emph{D}^{*}((H_{\bcP}^{\bS}(n,\bx))_{n=1}^{M}).
\end{equation} 
Hence the Theorem is proved. \qed \\

{\bf Bibliography.}

[BeCh]
 Beck, J., Chen, W.~W.~L., 
 Irregularities of Distribution,  
 Cambridge Univ. Press,  Cambridge, 1987.

 [Bi]
Bilyk, D.,  On Roth's orthogonal function method in discrepancy theory. Unif. Distrib. Theory, 6 (2011), no. 1, 143-184.

[Fa] Faure, H., Discr\'{e}pances de suites associ\'{e}es a un syst\`{e}me de num\'{e}ration (en dimension
un). Bull. Soc. Math. France 109: 143-182, 1981. (French)

[FaCh]
Faure, H., Chaix, H., Minoration de discr\'{e}pance en dimension deux, Acta Arith. 76 (1996), no. 2, 149-164. 

[FKP]  Faure, H., Kritzer, P., Pillichshammer F., 
From van der Corput to modern constructions of sequences for quasi-Monte Carlo rules,
 arXiv 1506.03764.

[He]
Hellekalek, P., Regularities in the distribution of special sequences, 
J. Number Theory, 18 (1984), no. 1, 41-55.

[Ha] Halton, J.H., On the efficiency of certain quasi-random sequences of points inevaluating
multi-dimensional integrals, Numer. Math., 2  (1960) 84-90.

[Le1]  Levin, M.B.,  On the lower bound in the  lattice point remainder problem  for a parallelepiped, to appear in Discrete \& Computational Geometry, arXiv: 1307.2080.

[Le2]  Levin, M.B.,  On the lower bound of the discrepancy  of Halton's sequences: I, 
  arXiv:1412.8705

[Le3] Levin, M.B.,  On the lower bound of the discrepancy  of $(t,s)$ sequences: II, arXiv: 1505.04975v2.	

[Ni]  Niederreiter, H., Random Number Generation and Quasi-Monte Carlo Methods, in: CBMS-NSF Regional Conference Series
in Applied Mathematics, vol. 63, SIAM, 1992.

[VC] van der Corput, J.G., Verteilungsfunktionen I-II. Proc. Akad. Amsterdam, 38 (1935), 813-
821, 1058-1066.
\end{document}